\documentclass[12pt,a4paper]{article}

\usepackage[latin1]{inputenc} 
\usepackage{amsmath, amssymb, amscd} 
\usepackage{theorem}

\newtheorem{teo}{Theorem}[section] 
\newtheorem{lemma}[teo]{Lemma} 
\newtheorem{cor}[teo]{Corollary} 
 
\theorembodyfont{\normalfont} 
\newtheorem{defi}[teo]{Definition} 
\newtheorem{rem}[teo]{Remark}

\newtheorem{exa}[teo]{Example}

\newcommand{\re}{\mathbb{R}}

\newcommand{\proj}{\mathbb{P}} 
\newcommand{\cpx}{\mathbb{C}}

\newcommand{\raw}{\rightarrow} 
 
\newcommand{\olo}{\mathcal{O}} 
 
\newcommand{\tr}{\textrm{tr}}

\newcommand{\diag}{\textrm{diag}}

\newcommand{\Bl}{\textrm{Bl}} 
\newcommand{\Spec}{\textrm{Spec}}
 
\newcommand{\barr}{\overline}

\newcommand{\del}{\partial} 
\newcommand{\delbar}{\barr{\del}}

\numberwithin{equation}{section}

\newcommand{\om}{\omega} 
 
\newcommand{\Lie}{\mathfrak{Lie}} 
 
\newcommand{\Aut}{\operatorname{Aut}}

\newcommand{\dimo}[1][]         {\noindent\textbf{Proof#1}. }

\newcommand{\restr}[1]          {\phantom{}_{\text{\raisebox{.4ex}{$|$}}#1}}

\newcommand{\fine}    {\begin{flushright} 
             $\Box$
             \end{flushright}}

\begin{document}
\title{Unstable blowups}

\author{Jacopo Stoppa} 
 
\date{} 
 
\maketitle

\begin{abstract} 
  \noindent 
Let $(X,L)$ be a polarised manifold. We show that K-stability and asymptotic Chow stability of the blowup of $X$ along a 0-dimensional cycle are closely related to Chow stability of the cycle itself, for polarisations making the exceptional divisors small. This can be used to give (almost) a converse to the results of Arezzo and Pacard in \cite{arezzoI}, \cite{arezzoII} and to give new examples of K\"ahler classes with no constant scalar curvature representatives.    
\end{abstract} 
\section{Introduction}\label{intro}
The theme of this paper is to construct and study \textit{test configurations} (i.e. particular degenerations) for blowups of a polarised manifold. More precisely we compute the \textit{Donaldson-Futaki invariant} of these configurations (Theorem \ref{main_teo}). These concepts are recalled in sections \ref{config} and \ref{futaki} respectively, but in essence the Donaldson-Futaki invariant is a rational number attached to our degeneration which morally plays the role of the Hilbert-Mumford weight in Geometric Invariant Theory. The main application is to nonexistence of constant scalar curvature K\"ahler (cscK for brevity) metrics contained in particular K\"ahler classes (\ref{arezzo_converse}). The search for such manifolds has played an important role in the development of the theory (see \cite{ross-thomas} for more details) and the method presented here yields infinitely many (which we can take to be rational, see \ref{exa_proj}).

A more algebraic application is to asymptotically Chow unstable polarisations, see \ref{chow_main_teo}.

Test configurations for a polarised manifold $(M,L)$ together with generalised Futaki invariants were introduced by Donaldson in \cite{donaldson_configurations} following work of Tian \cite{tian}. Let $S(\om)$ denote the scalar curvature of a K\"ahler metric $\om \in c_1(L)$ with average $\widehat{S} = 2 \pi n \frac{c_1(M)\cup c_1(L)^{n-1}}{c_1(L)^n}$. In \cite{donaldson_calabi} Donaldson proves that the \textit{Calabi functional}
\[Ca(\om)=\int_M (S(\om)-\widehat{S})^2 \,\frac{\om^n}{n!}\]
is bounded below in all of $c_1(L)$ by the negative of the generalised Futaki invariant $F(\mathcal{M})$ (divided by a positive term) for \textit{any} test configuration $(\mathcal{M},\mathcal{L})$ relative to $(M,L)$. Thus if $c_1(L)$ contains a cscK metric, $M$ must be \textit{K-semistable} with respect to $L$, meaning precisely $F(\mathcal{M}) \geq 0$ for all $(\mathcal{M}, \mathcal{L})$. Moreover it is expected that $F(\mathcal{M}) = 0$ if and only if $\mathcal{M}$ is isomorphic to the product $M \times \cpx$, i.e. a cscK manifold should be $\textit{K-polystable}$. That the converse holds as well is the content of the conjectural Hitchin-Kobayashi correspondence for manifolds (or Yau-Tian-Donaldson conjecture).
  
Now to our case. Let $X$ be a compact connected complex manifold, $\dim(X) = n$. We assume that $X$ is polarised, that is we fix an ample line bundle $L \raw X$. Then for any subscheme $Z \subset X$, the blowup $\Bl_Z X$ with exceptional divisor $E$ in endowed with the line bundle $L^{\gamma} - E$ which is ample for all large enough positive integers $\gamma$. Suppose now that $\Aut(X)$ contains a nontrivial compact connected subgroup, so that there is a 1-parameter subgroup (1-PS) of automorphisms of $X$ ( i.e. a group homomorphism $\alpha: \cpx^* \hookrightarrow \Aut(X)$). Then taking the limit of the action of $\alpha$ on $X$ as $t \raw 0$ induces in a natural way a test configuration $\mathcal{X}$ for $(\Bl_Z X, L^{\gamma}-E)$ intuitively by making the components of $E$ move around and possibly collide. In general one would like to make this rigorous and to understand this configuration as much as is needed to compute the first terms of the asymptotic expansion of its Futaki invariant $F(\mathcal{X})$ as $\gamma \raw \infty$. 
\begin{rem}
It is important to note that the we cannot expect the central fibre $\mathcal{X}_0$ to be the blowup of $X$ along the limit $Z_0$ of $Z$ as $t \raw 0$ in the relevant Hilbert scheme. This is explained at the end of section \ref{config}. In general we can only say that $\Bl_{Z_0}X$ is an irreducible component of $\mathcal{X}_0$. 
\end{rem}
Our main result is that we can carry out this program completely when $Z$ is a 0-dimensional cycle, say 
\[ Z = \sum_i a_i p_i\, (p_i \in X, a_i > 0).\]
\begin{rem}\label{multiplicity}
We emphasise that by blowing up $a_i p_i$ we mean blowing up the ideal $\mathcal{I}^{a_i}_{p_i}$, so that the \textit{algebraic} multiplicity of $p_i$ is the length of $\olo_{a_i p_i}$ :
\[\textrm{len}(\olo_{a_i p_i}) = {n + a_i-1\choose a_i -1}\]
while $a_i$ is the multiplicity of $p_i$ in the cycle. Note that by \cite{hart}, II Exercise 7.11, there is an isomorphism of polarised schemes
\[(\Bl_Z X, L^{\gamma}-E) \cong (\Bl_{\{p_i\}}, L^{\gamma} - \sum_i a_i E_i)\]
where $E_i$ is the component of the exceptional divisor over $p_i$. 
\end{rem}
In particular, we will relate $F(\mathcal{X})$ to two fundamental weights associated to the action of $\alpha$ on $X$: on the one hand, the classical Futaki invariant $F(X)$ for the holomorphic vector field which generates $\alpha$, on the other, the natural GIT weight for the action of $\alpha$ on the Chow variety of 0-cycles of total multiplicity $m = \sum_i a^{n-1}_i$ on $(X,L^{\gamma})$ (i.e. the symmetric product $X^{(m)}$ with polarisation induced by $(L^{\gamma})^{\boxtimes m}$). We denote this weight by $\mathcal{CH}(\sum_i a^{n-1}_i p_i,\alpha)$. The few GIT notions we need (including Chow stability for 0-cycles) are recalled in section \ref{chow_theory}.  
The following results are proved in section \ref{futaki}.
\begin{teo}\label{main_teo}
\[F(\mathcal{X}) = F(X)\,\gamma^{n} - \mathcal{CH}(\sum_i a^{n-1}_i p_i, \alpha)\frac{\gamma}{2(n-2)!} + O(1).\]
As a consequence if $(X, L)$ is K-polystable, the blowup of $X$ along a Chow-unstable 0-cycle is K-unstable when the exceptional divisors are small enough (i.e. $\gamma \gg 0$). 
\end{teo}
\begin{cor}\label{chow_main_teo}
If $(X,L)$ is asymptotically Chow polystable its blowup along a Chow unstable 0-cycle is asymptotically Chow unstable, when the exceptional divisors are small enough. 
\end{cor}
\begin{rem}\label{balancedness}
Thus when we blow up the cycle $\sum_i a_i p_i$ it is the GIT stability of the cycle $\sum_i a^{n-1}_i p_i$ that naturally shows up in $F(\mathcal{X})$. To interpret this difference note that the volume of the weighted exceptional divisor $a_iE_i$ over $p_i$ with respect to $L^{\gamma}-E$ is $\gamma^{1-n}a^{n-1}_i$ (up to a dimensional constant). Then an example in \cite{richard_notes} 27--28 illustrating Donaldson's theory of balanced metrics from \cite{donaldson_scalarI} suggests that when blowing up we are perturbing the centre of mass of $X$ in $\mathfrak{su}(H^0_X(L^{\gamma r})^*)$ (for $r, \gamma \gg 0$) by attaching a small weight proportional to $\gamma^{1-n}a^{n-1}_i$ over the point $p_i$.  
\end{rem}
We must mention at this point that our interest in this topic came from trying to find an algebro-geometric counterpart to the results of Arezzo and Pacard on blowing up and desingularizing cscK metrics contained in \cite{arezzoI}, \cite{arezzoII}. As blowing up is such a fundamental tool it is not surprising that the Arezzo-Pacard theorem plays a key role in the proofs of many recent results in the field. Among them we recall a theorem of Shu \cite{shu} stating that any compact complex surface with $b_1 = 0$ (except the blowup of $\proj^2$ in 1 or 2 points) is deformation equivalent to one bearing a cscK metric and the construction of Einstein metrics which are conformally K\"ahler on $\proj^2$ blown up in 2 points due to Chen-LeBrun-Weber \cite{chen}. By analogy it seems that understanding the behaviour of the Donaldson-Futaki invariant under blowing up may be useful for the algebraic side of the story.    

Working in the K\"ahler setting, Arezzo and Pacard prove that if $\omega$ is cscK on $X$, and if the cycle $\sum_i a_i p_i$ ($a_i \in \re^{+}$) satisfies the three conditions recalled below, then for all small positive $\epsilon$ the K\"ahler class on the blowup of $X$ at $\{p_i\}$ given by
\[\pi^*[\omega]-\epsilon\left(\sum a_i [E_i]\right)\] 
contains a cscK form $\omega_{\epsilon}$. Moreover $\omega_{\epsilon}$ converges to $\omega$ in the $C^{\infty}$ sense over $X-\sum_i p_i$.
These conditions are expressed in terms of a moment map
\[\mu: X \raw \mathfrak{ham}^*(X,J,\omega)\]
for the action of the Hamiltonian isometries of the K\"ahler manifold $X$, namely
\begin{enumerate}
\item $\{\mu(p_i)\}$ span $\mathfrak{ham}^*(X,J,\omega)$\label{arezzo_conditionI};
\item no nonzero element of $\mathfrak{ham}(X, J, \omega)$ vanishes at all the $\{p_i\}$\label{arezzo_conditionII}; 
\item $\sum_i a^{n-1}_i\mu(p_i) = 0$. \label{arezzo_conditionIII}
\end{enumerate}
\begin{rem}
Condition \ref{arezzo_conditionII} is only needed to get rid of residual automorphisms on the blowup, and without it the theorem continues to hold in a slightly modified form. It is also expected that \ref{arezzo_conditionI} is not necessary, but rather a side effect of the analytic argument used in the proof (see \cite{arezzoIII} for more details). Thus \ref{arezzo_conditionIII} (which we may call a balanced- or stability- condition) seems to be heart of the matter. This fits in well with the results of this paper. 
\end{rem}
We are now going to recast the Arezzo-Pacard theorem in terms of Chow stability. This bears on the projective case when $\omega = c_1(L)$ and the $a_i$ are positive integers. The Kempf-Ness theorem shows that the cycle $\sum_i a^{n-1}_i p_i$ can be modified by elements of $\Aut(X)$ to satisfy conditon \ref{arezzo_conditionIII} if and only if it is Chow polystable with respect to the action of $\Aut(X)$ (for a moment we refer to this new cycle as the balanced image).
\begin{rem}\label{kempf_ness}
One must be careful in applying the Kempf-Ness theorem here since the relevant symmetric product $X^{(m)}$ is a singular variety when $n > 1$, $m > 1$. However (see section \ref{chow_theory}) $X^{(m)}$ is defined as the geometric quotient $X^m /\Sigma_m$. Also the polarisation on $X^{(m)}$ is induced by $(L^{\gamma})^{\boxtimes m}$ on $X^m$. Since all the points of $X^m$ are stable under the action of $\Sigma_m$, we can lift a cycle in $X^{(m)}$ to the product $X^m$, apply Kempf-Ness there, and project back. 
\end{rem}
In the polystable case we might leave $\sum_i a^{n-1}_i p_i$ fixed and pullback $\omega$ by some automorphism instead so that \ref{arezzo_conditionIII} holds. Note that we can restrict to the action of the connected component of the identity $\Aut^0(X)$ in the Kempf-Ness theorem, and this preserves the K\"ahler class. This gives the following version of the theorem.
\begin{teo}
(Arezzo-Pacard, projective case) Suppose $\omega$ is cscK, the cycle $\sum_i a^{n-1}_ip_i$ is Chow polystable and conditions \ref{arezzo_conditionI}, \ref{arezzo_conditionII} hold on its balanced image. Then for all rational $\epsilon > 0$ small enough, the class
\[\pi^* [\omega] - \epsilon\left(\sum_i a_i E_i\right)\]
contains a cscK metric $\omega_{\epsilon}$. Moreover there exists $\phi \in \Aut^0(X)$ such that $\omega_{\epsilon}$ converges to $\phi^*\omega$ in the $C^{\infty}$ sense over $X - \sum_i p_i$.
\end{teo}  
Now observe that when $\omega$ is cscK $F$ vanishes identically on $\Lie\Aut(X)$. By the Hilbert-Mumford criterion we also know that $\sum_i a^{n-1}_i p_i$ is semistable if and only if $\mathcal{CH}(\sum_i a^{n-1}_i p_i,\alpha) \leq 0$ for all 1-PS $\alpha \hookrightarrow \Aut(X)$, so Theorem \ref{main_teo} almost implies that condition \ref{arezzo_conditionIII} is necessary. A small discrepancy comes from the fact that \textit{unstable} means \textit{not semistable} which is stronger than \textit{not polystable}. More precisely, our asymptotic expansion \ref{main_teo} together with Donaldson's lower bound on the Calabi functional immediately give 
\begin{teo}\label{arezzo_converse}
If $\omega$ is cscK and the cycle $\sum_i a^{n-1}p_i$ is Chow unstable then for all rational $\epsilon > 0$ small enough, the class
\[\pi^* [\omega] - \epsilon\left(\sum_i a_i E_i\right)\]
does not contain a cscK metric.
\end{teo}
\begin{exa}\label{exa_proj}
\textit{Projective space.} There is a very nice geometric criterion for Chow stability of points in $\proj^n$ (see for example \cite{mukai} 231--235). A cycle $Z = \sum_i m_i p_i$ is Chow \textit{unstable} if and only if for some proper subspace $V \subset \proj^n$ one has
\[\frac{|V \cap Z|}{\textrm{dim}(V)+1} > \frac{\sum_i m_i}{n+1}.\]
So already in the case of $\proj^2$ \ref{main_teo} gives infinitely many new examples of K-unstable classes: it suffices that more than $2/3$ of the points (counted with multiplicities) are aligned to get a K-unstable blowup. This gives a full generalisation of Example 5.30 in \cite{ross-thomasII} where it is shown that when $m_1 \gg m_j\, (j > 2)$ the blowup is K-unstable (actually \textit{slope unstable} with respect to $E_1$).    
\end{exa}
The following examples are straightforward applications of the theory in section \ref{chow_theory}.
\begin{exa}
\textit{Products.} Consider the product $X \times Y$ of $(X,L_X)$ and $(Y, L_Y)$ polarised by $L_X \boxtimes L_Y$. Then for $\alpha: \cpx^* \hookrightarrow \Aut(X)$,  \[\mathcal{CH}(\sum_i a_i p_i, \alpha \times 1) = \mathcal{CH}(\sum_i a_i \pi_X(p_i), \alpha).\]
Thus we may apply the above geometric criterion to the fibres of a product $\proj^n \times Y$. For example if $Y$ is K-polystable the blowup of $\proj^n \times Y$ along an unstable cycle supported at a single $\proj^n$-fibre will be K-unstable. A special case is the product $\proj^n \times \proj^m$ polarised by $\olo(1)\boxtimes\olo(1)$. A $0$-cycle will be unstable whenever its projection to one of the two factors is, e.g. in the case of 3 distinct point, when 2 of them lie on a vertical or horizontal fibre. This gives more examples of unstable blowups.  
\end{exa}
\begin{exa}\textit{$\proj^1$ bundles.} Similarly we can consider the projective completion $X$ of some line bundle $L$ over a polarised manifold. In this case the so-called momentum construction yields many examples of cscK metrics, see \cite{hwang}. Any polarisation $\mathcal{L}$ on $X$ restricts to $\olo_{\proj^1}(k)$ on all the fibres for some $k$. There is a natural $\cpx^*$-action on $L$ given by complex multiplication on the fibres, and this extends to a $\cpx^*$-action $\alpha$ on $X$ so that points lying on the zero (resp. infinity) section $X_0$ ($X_{\infty}$) are fixed, but with weight $k$ (resp. $-k$) on the line above them. By acting with $\alpha^{-1}$ instead if necessary, we conclude that whenever more than half the points lie on $X_0$ or $X_{\infty}$, the corresponding $0$-cycle is Chow unstable (i.e. its Chow weight is $ > c k$ for some positive constant $c$). We see that the conclusion is really independent of $k$ and so the blowup along such a cycle will be K- and asymptotically Chow- unstable for \textit{any} polarisation on $X$ making the base cscK.
\end{exa}
\textbf{Notation.} We will often suppress pullback maps and use the same letter to denote a divisor and the associated line bundle. Consequently we mix additive and multiplicative notation as necessary.\\
\\ \textbf{Acknowledgements}. This work has been made possible by the generosity and patience of my supervisor R. Thomas. The idea of this project is due to him and S. K. Donaldson who I also thank for his kind interest. I would like to thank J. Ross, G. Székelyhidi and the Referee for many important suggestions. Finally I am grateful to the geometry groups at Imperial and Pavia and in particular to G. P. Pirola and A. Ghigi for their support and for teaching me so much.  
\section{Test configurations coming from automorphisms of the base}\label{config}
Let $Z = \sum_i a_i p_i$ be some 0-dimensional cycle on $X$, that is the closed subscheme supported at the points $\{p_i\}$ with nonreduced structure $\prod_i\mathcal{I}^{a_i}_{p_i}$. In this section we construct a \textit{test configuration} for $(\Bl_Z X, L^{\gamma}-E)$ naturally associated to a 1-PS $\alpha \hookrightarrow \Aut(X)$. A test configuration for a polarised manifold $(M, L)$ is given by a polarised flat family $(\mathcal{M}, \mathcal{L})\raw \cpx$ endowed with an $\mathcal{L}$-linearised $\cpx^*$-action covering the usual action of $\cpx^*$ on $\cpx$, and such that for any $t \neq 0$, $(\mathcal{M}_t, \mathcal{L}\restr{t}) \cong (M, L^s)$ for some exponent $s$ (called the exponent of the test configuration). In our case this is the natural flat family induced by the $\cpx^*$-action on $X$ and so on cycles and exceptional divisors. It will be useful to introduce the flat family of closed subschemes of $X$ given by $\{Z_t = \alpha(t)Z$, $t \in \cpx^*\}$.
\begin{lemma}\label{flat_family}
There is a flat family $p: \mathcal{X} \raw \cpx$ such that $\mathcal{X}_t \cong \emph{Bl}_{Z_t} X$ for all $t \in \cpx - \{0\}$. This is endowed with an induced action of $\alpha$ covering the usual $\cpx^*$ action on $\cpx$. 
\end{lemma}
\dimo We see $Z$ as a point of the Hilbert scheme of closed subschemes of $X$ with \textit{constant} Hilbert polynomial equal to the length of $\olo_Z$ as a module over itself. From the general theory we know that the flat family $(Z_t,t) \subset X \times (\cpx-\{0\})$ has a unique flat closure i.e. there exists a unique closed subscheme $Y \subset X \times \cpx$, flat over $\cpx$, with fibres $Y_t = Z_t$ for all $t \in \cpx-\{0\}$. As total scheme of our test configuration we take $\mathcal{X} = \Bl_Y(X\times\cpx)$, with projection $p: \mathcal{X} \raw \cpx$ given by the composition of the blowup map $\pi: \mathcal{X} \raw X \times \cpx$ with the projection onto the second factor. By \cite{hart}, II Proposition 7.16 $\mathcal{X}$ is reduced, irreducible, and $p$ is a dominant (in fact surjective) morphism. The base $\cpx$ is regular, $1$-dimensional and dominated by every irreducible component of $\mathcal{X}$ so by \cite{hart}, III Proposition 9.7 we see that $p$ is a flat morphism. Since $Y$ is preserved by $\alpha$ there is an induced action of $\alpha$ on $\mathcal{X}$ covering the usual action of $\cpx^*$ on $\cpx$. \fine
Since $\mathcal{X} = \textrm{Proj}(\bigoplus_r \mathcal{I}^r_Y)$, it is naturally endowed with an invertible sheaf $\olo(1)$. Let $p_X: \mathcal{X} \raw X$ be the composition $\mathcal{X} \raw X\times \cpx \raw X$. Define a line bundle on $\mathcal{X}$ by $\mathcal{L} = p_X^*L^{\gamma}\otimes\olo(1)$. A slight modification of the above argument then proves
\begin{lemma}\label{line_bundle}
For all large $\gamma$, $\mathcal{L}$ is a $\alpha$-linearised ample line bundle $\mathcal{L}$ on $\mathcal{X}$ such that $\mathcal{L}\restr{t} \cong L^{\gamma} - E$ for $t \neq 0$.
\end{lemma}
Next we need to study the central fibre $\mathcal{X}_0$ of our test configuration. It will be useful to write $Z_0$ for the limit of $Z_t$ in the Hilbert scheme as $t \raw 0$ as in the proof of \ref{flat_family}, and to define
\[\widehat{X} = \Bl_{Z_0}X\]
and 
\[E_0 = \text{the exceptional divisor of  } \widehat{X} \raw X.\]
We claim that there is a closed immersion $\widehat{X} \hookrightarrow \mathcal{X}_0$. To see this, consider the closed immersion $i: X \cong X \times \{0\} \hookrightarrow X \times \cpx$. By flatness we see that
\begin{lemma}\label{inverse_image}
The inverse image ideal sheaf $i^{-1}\mathcal{I}_Y\cdot\olo_X $ is $\mathcal{I}_{Z_0}$.
\end{lemma}
\begin{lemma}
The inclusion $i: X \times {0} \hookrightarrow X\times\cpx$ induces a closed immersion
\[\widehat{i}: \widehat{X} \hookrightarrow \mathcal{X}_0.\]
\end{lemma}
\dimo
By \cite{hart},II Corollary 7.15, there is an induced closed immersion
\[\widehat{i}: \Bl_{i^{-1}\mathcal{I}_Y\cdot\olo_X} X \hookrightarrow \Bl_{Y}(X\times\cpx).\]
Since the image lies in $\mathcal{X}_0$ we conclude by \ref{inverse_image}.
\fine
Next we define a closed subscheme $P$ of $\mathcal{X}_0$ by 
\[P = (\mathcal{X}_0 - \widehat{X})^{-}.\]
In general $\widehat{i}$ is not an isomorphism as the exceptional set $P$ may well be a component of $\mathcal{X}_0$. To motivate this recall that $Z_0$ is the central fibre of the flat family $\{Z_t, t \in \cpx\}$. But a family of thickenings $\{r Z_t, t \in \cpx\}$ will \textit{not} be flat in general. The generic fibre of $\mathcal{X}$ is $\mathcal{X}_t = \textrm{Proj}(\bigoplus_{r} \mathcal{I}_{r Z_t})$ for $t \neq 0$. As $\mathcal{X}$ itself is flat we see that $\mathcal{X}_0$ cannot in general be $\textrm{Proj}(\bigoplus_{r} \mathcal{I}_{r Z_0})$. This means there is some extra closed subscheme $P$ inside $\mathcal{X}_0$, which we will have to take into account in section \ref{futaki} when computing the Futaki invariant.
\begin{exa}
As an affine example consider the ideal \[\mathcal{I}_Y  = (x(x-t), xy, y(y-t)) \subset \cpx[x,y,t]\]
which describes 3 points colliding along 2 orthogonal directions in $\mathbb{A}^2$. One can show that $\mathcal{X}$ is the closed subscheme of $\Spec\,\cpx[x,y,t]\times\textrm{Proj}\,\cpx[\xi_0,\xi_1, \xi_2]$ defined by
\[\mathcal{I}_{\mathcal{X}} = ((x-t)\xi_1 - y\xi_0, (y-t)\xi_1 - x\xi_2)\]
with central fibre (sitting inside $\Spec\,\cpx[x,y]\times\textrm{Proj}\,\cpx[\xi_0,\xi_1, \xi_2]$)
\[\mathcal{I}_{\mathcal{X}_0} = (x\xi_1 - y\xi_0, y\xi_1 - x\xi_2).\]
In this case $\widehat{X}$ is the closed subscheme of $\mathcal{X}_0$ given by
\[\mathcal{I}_{\widehat{X}} = (x\xi_1-y\xi_0, y\xi_1 - x\xi_2, \xi^2_1-\xi_1\xi_2)\]
that is a whole component of $\mathcal{X}_0$, while the exceptional component $P$ is $\textrm{Proj}\,\cpx[\xi_0,\xi_1, \xi_2]\cong\proj^2$. 
\end{exa}
In any case the restriction of $\mathcal{L}_0$ to $\widehat{X}$ is the expected one. 
\begin{lemma}\label{central_line}
\[\mathcal{L}_0\restr{\widehat{X}} = L^{\gamma} - E_0.\]
\end{lemma}
The easy proof is left to the reader.
\section{Chow stability for 0-dimensional cycles}\label{chow_theory}
In this section we recall the few GIT notions we need in the case of the Chow variety of points on $X$. For much more on this see \cite{mum} Chap. 3. Let the symmetric group $\Sigma_d$ on $d$ letters acts on the $d$-fold product $X^d$. The symmetric product $X^{(d)} = X^d/\Sigma_d$ is a projective variety. The points of $X^{(d)}$ are the orbits of the $d$-tuples of points of $X$ under permutation and so can be identified with effective $0$-cycles $\sum n_i [x_i]$ with $x_i \in X$, $n_i > 0 $ and $\sum n_i = d$. This shows that $X^{(d)}$ is actually the Chow variety of length $d$ 0-cycles on $X$. The construction of an $\Aut(X)-$linearised ample line on $X^{(d)}$ can be made very explicit as follows. Let $V = H^0(X,L^{\gamma})^*$ for some large ${\gamma}$ and embed $X \hookrightarrow \proj(V)$. Denote by $\proj(V^*)$ the projective space of hyperplanes in $\proj(V)$, and by $Div^d(\proj(V^*))$ the projective space of effective divisors of degree $d$ in $\proj(V^*)$. For any $p \in \proj(V)$ consider the hyperplane in $\proj(V^*)$ given by
\[H_p := \{l \in \proj(V^*) : p \in l\}.\]
Define a morphism $ch: (\proj(V))^d \raw Div^d(\proj(V^*))$ by 
\[ ch(x_1, ..., x_d) := \sum_i H_{x_i}.\]
This is the \textit{Chow form} of $\{x_1,...,x_d\}$: the
divisor of hyperplanes whose intersection with $\{x_1,...,x_d\}$ is
nonempty. As for $X^d$ we have the composition
\[ch: X^d \hookrightarrow \proj(V)^d \hookrightarrow Div^d(\proj(V^*))\]
induced by the product line $(L^{\gamma})^{\boxtimes d}$. Now $ch$ is $\Sigma_d$ equivariant and by the universal property of the geometric quotient factors through $\proj(V)^{(d)}$ defining a morphism 
\[ch: X^{(d)} \raw Div^d(\proj(V^*)).\]
One can check that $ch$ defines an isomorphism on its image. We can thus identify $X^{(d)}$ with its image $ch(X^{(d)})$ in $Div^d(\proj(V^*))$. By the usual identifications 
\[\label{divisors} Div^d(\proj(V^*)) \cong \proj(H^0(\proj(V^*), \olo(d))) \cong \proj(S^d V) \]
we see $ch$ as a map with values in $\proj(S^d V)$:
\[ch: X^{(d)} \hookrightarrow \proj(S^d V).\]
Under these identifications then $ch$ is the map
\begin{equation}\label{CH_weight}
X^{(d)} \ni \{[x_1],...[x_d]\} \mapsto [x_1\cdot...\cdot x_d]
\end{equation}
defined via the embedding $X \hookrightarrow \proj(V)$. 
\begin{rem}
It is important to emphasise that $ch$ is given by the descent of $(L^{\gamma})^{\boxtimes d}$ under the action of $\Sigma_d$. This means that the Chow line $\olo_{Div^d(\proj(V^*))}(1)\restr{X^{(d)}}$ pulls back to $(L^{\gamma})^{\boxtimes d}$ under the quotient map. This holds because by \ref{CH_weight} $\olo_{Div^d(\proj(V^*))}(1)$ pulls back to the line $\olo_{\proj(V)}(1)^{\boxtimes d}$ on $\proj(V)^d$ under $ch$ and this in turn pulls back to $(L^{\gamma})^{\boxtimes d}$.
\end{rem}
Now we assume that $\alpha \hookrightarrow \Aut(X)$ acts through a 1-PS $\alpha \hookrightarrow \textrm{Sl}(V)$ and we come to the definition of the GIT weight for the action of $\alpha$ on $X$. Recall that we can find a basis of eigenvectors so that $\alpha(t)$ acts as $\diag(t^{\lambda_0}, ..., t^{\lambda_N})$, where $\dim(V) = N+1$. In these projective coordinates on $\proj(V)$, writing $X \ni x = [v_0:...:v_n]$, we define the Mumford weight as 
\[\lambda(x, \alpha) := \textrm{min}\{\lambda_i:v_i \neq 0\}.\] 
By definition of $V$ this is the weight of the induced action on the line $(L^{\gamma})^*$ over the limit $x_0$ of $\alpha(t)x$ as $t \raw 0$ in the $\cpx^*$-action.

In the same way we can define the Mumford weight for the induced action of $\alpha$ on $X^{(d)}$. In fact $\alpha \hookrightarrow \textrm{Sl}(V)$ naturally induces a 1-PS $\alpha \hookrightarrow \textrm{Sl}(S^d V)$. Thus the embedding $ch$ described above gives a natural linearisation for this action. We write $\mathcal{CH}$ for this Mumford weight. Then by \eqref{CH_weight} we immediately obtain the relation
\[\mathcal{CH}(\sum_i m_i x_i, \alpha) = \sum_i m_i \lambda(x_i, \alpha)\]
for any $\sum_i m_i x_i \in X^{(d)}$. 
\begin{rem}\label{gamma_chow}
While the numerical value of the Mumford weight depends on the power $L^{\gamma}$ we take, the fact that a cycle $\sum_i m_i p_i$ is semistable is independent of ${\gamma}$. This is immediate from the definition of stability in terms of invariant sections.
\end{rem}

\section{The Donaldson-Futaki invariant}\label{futaki}
In this section we compute the Donaldson-Futaki invariant of the induced $\cpx^*$ action on the central fibre $\mathcal{X}_0$. But first we recall Donaldson's definition of the Futaki invariant of a $\cpx^*$-action on a variety (or scheme) $M$ of dimension $n$ endowed with a linearised ample line bundle $L$. We write $A_k$ for the infinitesimal generator of the induced $\cpx^*$-action on $H^0(M, L^k)$. 
\begin{rem}
The lifting of the $\cpx^*$-action to $L$ is \emph{not} unique, so $A_1$ is not well defined. However for any other lifting, there is $\lambda \in \mathbb{Z}$ such that
\[A'_1 = A_1 + \lambda I_1.\]
where $I_1$ denotes the identity matrix on $H^0(M,L)$. As a consequence, 
\[A'_k = A_k + k \lambda I_k\]
where $I_k$ denotes the identity matrix on $H^0(M,L^k)$. Using this transformation rule one can check that the Futaki invariant defined below is independent of the choice of lifting. 
\end{rem}
By Riemann-Roch and its equivariant version there are expansions
\[h^0(M,L^k) = c_0 k^n + c_1 k^{n-1} + O(k^{n-2}),\]
\[\tr(A_k) = b_0 k^{n+1} + b_1 k^{n} + O(k^{n-1})\]
valid for all large $k$. The Futaki invariant is defined as
\begin{equation}\label{futaki_def}
F = \frac{c_1 b_0}{c_0} - b_1.
\end{equation}
According to this general definition, in our case we need to compute $h^0(\mathcal{X}_0, \mathcal{L}^r_0)$ and the trace of the induced action on $H^0(\mathcal{X}_0, \mathcal{L}^r_0)$ for all large $r$. It is important to keep in mind that $\mathcal{L}_0$ depends in turn on the parameter $\gamma$; i.e. $\mathcal{L}_0$ comes from picking the line $L^{\gamma} - E$ on the generic fibre. 
\begin{rem}
We emphasise that with this choice of notation $\gamma^{-1}$ measures the volume of the exceptional divisors on the generic fibre, while $r$ is the scale parameter needed to compute the Futaki invariant. Also in the proof of the following lemma and in many other places below, we need the vanishing of higher cohomology groups. This will always hold for $r \geq r_0 = r_0(\gamma)$, however increasing $\gamma$ only makes $\mathcal{L}$ (and other related line bundles) more positive, so at least for $\gamma \geq \gamma_0$ we can assume that $r_0$ is a fixed constant.
\begin{lemma}\label{dimI}
\end{lemma}
\[h^0(\mathcal{X}_0, \mathcal{L}^r_0) = h^0(X, L^{\gamma r}) - \left(\sum_i{a^n_i}\right)\frac{r^{n}}{n!} -\left(\sum_i a^{n-1}_i\right)\frac{r^{n-1}}{2(n-2)!} + O(r^{n-2}).\]
\end{rem}
\dimo
By flatness $h^0(\mathcal{X}_0, \mathcal{L}^r_0) = h^0(\Bl_Z X, L^{\gamma r}-rE)$. To compute the coefficients use the asymptotic Riemann-Roch formula, keeping Remark \ref{multiplicity} in mind.
\fine
Now for the trace. To try to keep the notation light in what follows we will write $\tr(U)$ for the trace of the induced action on some vector space $U$. We start with the restriction $\cpx^*$-equivariant exact sequence 
\begin{equation}\label{restrictionI}
0 \raw H^0_P(\mathcal{I}^r_{E_0}\mathcal{L}^r_0\restr{P})\raw H^0_{\mathcal{X}_0}(\mathcal{L}^r_0) \raw H^0_{\widehat{X}}(L^{\gamma r} - r E_0)\raw 0
\end{equation}
which holds for large $r$. So we see 
\begin{equation}\label{big_trace}
\tr(H^0_{\mathcal{X}_0}(\mathcal{L}^r_0)) = \tr(H^0_{\widehat{X}}(L^{\gamma r} - r E_0)) + \tr(H^0_P(\mathcal{I}^r_{E_0}\mathcal{L}^r_0\restr{P})).
\end{equation}
To compute the first term we turn to the natural isomorphism (for $r \gg 0$)
\[H^0_{\widehat{X}}(L^{\gamma r} - r E_0) \cong H^0_X(\mathcal{I}^r_{Z_0}L^{\gamma r}).\]
The exact sheaf sequence on $X$
\begin{equation}
0 \raw \mathcal{I}^r_{Z_0}L^{\gamma r} \raw L^{\gamma r}\raw \olo_{r Z_0}\otimes_{\cpx} L^{\gamma r}\restr{Z_0}\raw 0
\end{equation}
is $\cpx^*$-equivariant and gives an exact sequence of sections (for large $r$):
\begin{equation}\label{restrictionII}
0 \raw H^0_{X}(\mathcal{I}^r_{Z_0}L^{\gamma r}) \raw H^0_X(L^{\gamma r})\raw \olo_{r Z_0}\otimes_{\cpx} L^{\gamma r}\restr{Z_0}\raw 0.\end{equation}
Here we used that $\olo_{r Z_0}$ is a skyscraper sheaf supported at $Z_0$. So we see
\begin{equation}
\tr(H^0_{\widehat{X}}(L^{\gamma r}-r E_0)) = \tr(H^0_X(L^{\gamma r})) - \sum_q \tr(\olo_{r Z_0,q} \otimes L^{\gamma r}\restr{q}) 
\end{equation}
where we are summing over the components of $(Z_0)_{\textrm{red}}$. Substituting in \eqref{big_trace} we get
\begin{equation}\label{trace}
\tr(H^0_{\mathcal{X}_0}(\mathcal{L}^r_0)) = \tr(H^0_X(L^{\gamma r})) - \sum_q \tr(\olo_{r Z_0,q} \otimes L^{\gamma r}\restr{q}) + \sum_q \tr(H^0_{P_q}(\mathcal{I}^r_{E_{0,q}}\mathcal{L}^r_0\restr{P_q})) 
\end{equation} 
where we write $P_q$ for the component of $P$ which projects to $q$ via $p_X$ and similarly for $E_{0,q}$. 
\begin{defi}
For any $q \in (Z_0)_{\textrm{red}}$, we denote by $\lambda(q)$ the weight of the induced $\cpx^*$-action on the line       $L\restr{q}$. Note that the weight of the induced action on $L^{m}\restr{q}$ is then $m \lambda(q)$ for any $m > 0$. As we already observed $\lambda(q)$ depends on the choice of a lifting of $\alpha$ to $L$, but remember that this choice will not affect the Futaki invariant. For any other lifting the new weights are 
\[\lambda'(q) = \lambda(q) + \lambda\]
for some $\lambda \in \mathbb{Z}$. These weights should not be confused with the relevant Chow weights, which require $\alpha$ to act through $\textrm{Sl}(H^0(X, L^m))$ (for the relevant power $m$) in their definition. This difference will turn out to be important for our purposes.
\end{defi}
The crucial step to get an asymptotic expansion for $F(\mathcal{X})$ is the following rough estimate, ignoring any term which is independent of $\gamma$.
\begin{lemma}\label{trI}
\[\emph{tr}(\olo_{r Z_0,q} \otimes L^{r \gamma}\restr{q}) = (r\gamma)\lambda(q)\emph{dim}(\olo_{r Z_0,q}) + O(\gamma^0 r^{n+1}).\]
\end{lemma}
\dimo
This is just the statement that the induced action on $\olo_{r Z_0, q}$ as a $\cpx$-vector space does not depend on the parameter $\gamma$.
\fine
There is a similar estimate for the action on the components $P_q$ of $P$.
\begin{lemma}\label{trII}
\[\emph{tr}(H^0_{P_q}(\mathcal{I}^r_{E_{0,q}}\mathcal{L}^r_0\restr{P_q})) = (r \gamma)\lambda(q)h^0_{P_q}(\mathcal{I}^r_{E_{0,q}}\mathcal{L}^r_0\restr{P_q}) + O(\gamma^0r^{n+1}).\]
\end{lemma}
\dimo
Note that
\[\mathcal{I}^r_{E_{0,q}}\mathcal{L}^r_0\restr{P_q} \cong L^{\gamma r}\restr{q} \otimes_{\cpx}\mathcal{I}^r_{E_{0,q}}\olo(r)\restr{P_q}\]
and that $L\restr{q}$ is the trivial line on $P_q$ (since it is pulled back from $Z_0$) acted on by $\cpx^*$ with weight $\lambda(q)$. \fine
For the following lemma we introduce the sets
\[A_q = \{p_i \in Z_{\textrm{red}}: \lim_{t\raw 0}\alpha(t)p_i = q\}.\]
\begin{lemma}\label{dimII}
\[\emph{dim}(\olo_{r Z_0,q}) = h^0_{P_q}(\mathcal{I}^r_{E_{0,q}}\mathcal{L}^r_0\restr{P_q}) + \left(\sum_{p_i \in A_q}a^{n}_i\right)\frac{r^n}{n!} + \left(\sum_{p_i \in A_q} a^{n-1}_i\right)\frac{r^{n-1}}{2(n-2)!} + O(r^{n-2}).\]
\end{lemma}
\dimo
This follows from local versions of \ref{dimI}, \eqref{restrictionI} and \eqref{restrictionII} around $q$ (in the analytic topology).
\fine
Putting these results together we can finally compute the trace on the central fibre.
\begin{lemma}\label{trIII}
\[\emph{tr}(H^0_{\mathcal{X}_0}(\mathcal{L}^r_0)) = \emph{tr}(H^0_X(L^{\gamma r}))\]
\[ -\left( \gamma \sum_q \lambda(q)\left(\sum_{p_i \in A_q}a^{n}_i\right)\frac{r^{n+1}}{n!}+\gamma \sum_q \lambda(q)\left(\sum_{p_i \in A_q} a^{n-1}_i\right)\frac{r^{n}}{2(n-2)!}\right)\]\[ + O(\gamma^0 r^{n+1}).\]
\end{lemma}
\dimo
Substitute the results of lemma \ref{trI} and lemma \ref{trII} into \eqref{trace} using \ref{dimII} to compute the missing dimension.
\fine
\begin{rem}\label{cancellation}
The reader may notice that this result depends crucially on the cancellation of the terms arising from $h^0_{P_q}(\mathcal{I}^r_{E_{0,q}}\mathcal{L}^r_0\restr{P_q})$ over the various points $q$. These are essentially (to higher order in $\gamma$) the terms encoding the singularities which form in $\mathcal{X}$ as $t \raw 0$, over which we have little control. A first reason why $F(\mathcal{X})$ should not ``see" these singularities is the argument with balanced metrics \cite{richard_notes} 27--28  already mentioned in \ref{balancedness} which we do not reproduce here. An alternative differential-geometric argument is sketched in \ref{diff_geo} below.
\end{rem}
The results obtained so far can be put in a form which makes applying definition \eqref{futaki_def} easier. Define co-efficients $b_i$, $c_i$ by $h^0_X(L^{\gamma r}) = c_0\gamma^n r^n + c_1\gamma^{n-1}r^{n-1} + O(r^{n-2})$, $\tr(H^0_X(L^{\gamma r})) = b_0\gamma^{n+1}r^{n+1} + b_1\gamma^{n}r^{n} + O(r^{n-1})$ . Similarly we define $b'_i$, $c'_i$ by $h^0_{\mathcal{X}_0}(\mathcal{L}^r_0) = c'_0(\gamma)r^{n} + c'_1(\gamma)r^{n-1} + O(r^{n-2})$, $\tr(H^0(\mathcal{X}_0, \mathcal{L}^r_0)) = b'_0(\gamma)r^{n+1} + b'_1(\gamma)r^{n} + O(r^{n-1})$.
\begin{cor}
\[b'_0 = b_0 \gamma^{n+1} - \sum_q \lambda(q)\left(\sum_{p_i \in A_q}a^{n}_i\right)\frac{\gamma}{n!}+O(1),\]
\[c'_1 = c_1 \gamma^{n-1} - \frac{1}{2(n-2)!}\sum_i{a^{n-1}_i},\]
\[c'_0 = c_0 \gamma^{n} - \frac{1}{n!}\sum_i{a^n_i},\]
\[b'_1 = b_1\gamma^{n} - \sum_q \lambda(q)\left(\sum_{p_i \in A_q} a^{n-1}_i\right)\frac{\gamma}{2(n-2)!}+O(1).\]
\end{cor}
\dimo
This is a restatement of \ref{dimI} and \ref{trIII}. 
\fine
With these preliminary computations in place we can now prove our main result.\\
\\
\textbf{Proof of \ref{main_teo}.}
By \eqref{futaki_def}
\[F(\mathcal{X}) = \frac{b'_0 c'_1}{c'_0}-b'_1 = \] \[=\left(\frac{b_0 c_1}{c_0} - b_1\right)\gamma^n \]\[+ \frac{1}{2(n-2)!}\left(\sum_q \lambda(q)\left(\sum_{p_i \in A_q} a^{n-1}_i\right)-\frac{b_0}{c_0}\sum_i{a^{n-1}_i}\right)\gamma + O(1)= \]
\[= F(X)\,\gamma^n + \frac{1}{2(n-2)!}\left(\sum_q \sum_{p_i \in A_q} a^{n-1}_i\left(\lambda(q) - \frac{b_0}{c_0}\right)\right)\gamma + O(1).\]
It remains to make the connection with Chow stability. Recall that to define the Chow weights with respect to the line $L^{\gamma}$ on $X$ we need $\alpha$ to act through $\textrm{Sl}(H^0_{X}(L^{\gamma})^*)$. Choosing any infinitesimal generator $A_{\gamma}$ for the action on $H^0_{X}(L^{\gamma})^*$ we need to solve for a correction parameter $\lambda_{\gamma}$
\[\tr(A_{\gamma}) + \gamma\lambda_{\gamma}\, h^0_{X}(L^{\gamma}) = 0\]
so we get
\[\lambda_{\gamma} = -\frac{\tr(A_{\gamma})}{\gamma h^0_{X}(L^{\gamma})} = -\frac{b_0}{c_0} + O(\gamma^{-1}).\]
Note that after pulling back the family $\mathcal{X}$ by a finite covering of $\cpx$ (i.e. $t\mapsto t^k$ for some $k$) we may assume $\lambda_{\gamma} \in \mathbb{Z}$. So by substituting $\lambda_{\gamma} + O(\gamma^{-1})$ for $-\frac{b_0}{c_0}$ the expansion above may be read as
\[F(\mathcal{X}) = F(X)\gamma^n + \frac{1}{2(n-2)!}\left(\sum_q \sum_{p_i \in A_q} a^{n-1}_i\lambda'(q)\right)\gamma + O(1)\]
where $\lambda'(q)$ are the new special linear weights 
\[\lambda'(q) = \lambda(q) + \lambda_{\gamma}.\]
By the discussion in section \ref{chow_theory} we see
\[\sum_q \sum_{p_i \in A_q} a^{n-1}_i\lambda'(q) = -\mathcal{CH}(\sum_i a^{n-1}_i p_i, \alpha)\]
where $\mathcal{CH}$ stands for the Chow weight relative to the polarisation on $X^{(d)}$ induced by $(L^{\gamma})^{\boxtimes d}$.
\fine
\begin{rem} In view of the proof we should write $\mathcal{CH}_{\gamma}$ in theorem \ref{main_teo}, but we drop the dependence on $\gamma$ motivated by Remark \ref{gamma_chow}. We should also write $\mathcal{X}_{\gamma}$ but this dependence is not really serious since we are only pulling back $\mathcal{X}$ by a finite covering.
\end{rem}
\textbf{Proof of \ref{chow_main_teo}.} By  \cite{ross-thomas} Theorem 3.9 and Theorem 4.33 we know that if $(X,L)$ is asymptotically Chow stable then it is K-semistable. In particular this implies $F(X) = 0$. Now we blowup along a Chow unstable cycle and apply \ref{main_teo} to get $F(\mathcal{X}) < 0$ for $\gamma \gg 0$, i.e. $(\Bl_{Z} X, L^{\gamma}-E)$ is K-unstable, and so in turn asymptotically Chow unstable.\fine
\begin{rem}\label{diff_geo}
We now give the argument promised in \ref{cancellation}. For this we need to recall the definition of the $K$-energy functional (due to Mabuchi \cite{mabuchi}). Consider the space of K\"ahler potentials with respect to a fixed K\"ahler form $\om$, $\mathcal{H} = \{\phi \in C^{\infty}(M, \re) : \om_{\phi} = \om + i\del\delbar\phi > 0\}$. $K$-energy is the unique functional $\mathcal{M}$ on $\mathcal{H}$ (up to an additive constant) such that the derivative $\frac{d}{dt}\mathcal{M}(t)$ along any path $\om_t = \om + i \del\delbar \phi_t$ is given by $\int_M \left(\frac{d}{dt}\phi_t\right) (S(\om_t)-\widehat S)\,\om^n_t$. Now by the general theory (markedly the moment map picture in \cite{don_fields}) one expects that the Donaldson-Futaki invariant for a test configuration $\mathcal{X}$ can be computed via any path of metrics $\om_t$, $t \in (0,1]$ which is ``adapted" to $\mathcal{X}$ (we will not try to make this precise here). The prediction is then $\lim_{t\raw 0}t\frac{d}{dt}\mathcal{M}(t) = F(\mathcal{X})$. Let us spell out what this means in our case, starting with a K\"ahler form $\om$ on the base $X$ and a holomorphic vector field $\alpha$ with flow $\alpha(t)$. Choose small enough \emph{disjoint} coordinate balls $B_i(2\epsilon)$ around each $p_i$. At a fixed time $t$, we condider the metric $\alpha(t)^*\om$ on $X$; remove a ball $B_i (\epsilon)$ and glue in a metric on a large open neighborhood of the zero section of $\olo(-1)\raw \proj^{n-1}$, such that the volume of the zero section is $\epsilon^{n-1}a^{n-1}_i$ (up to a constant). This will require deforming $\alpha(t)^*\om$, but we can leave it \emph{unchanged} outside $B_i(2\epsilon)$. This construction can be made to yield a sequence of smooth metrics on $\Bl_{\{p_i\}} X$ with the required cohomology class. Then taking the limit as $t \raw 0$ morally yields $F(\mathcal{X})$. On the other hand, it should be clear by construction and the definition of $\frac{d}{dt}\mathcal{M}(t)$ that this limit only depends on the action of $\alpha(t)$ in a neighborhood of each $p_i$ for small $t$; that is, the result does \emph{not} depend on mutual interaction of the points $\{p_i\}$ we blowup. While this conjectural picture could be overly difficult to make precise, it gives some geometric meaning to the cancellation of the $h^0_{P_q}(\mathcal{I}^r_{E_{0,q}}\mathcal{L}^r_0\restr{P_q})$ terms. 
\end{rem}
\begin{rem}
A striking feature of \ref{main_teo} is that it naturally suggests that higher order contributions to $F(\mathcal{X})$ should arise from blowing up higher dimensional subschemes.  
\end{rem}
\begin{rem}
Note that \ref{main_teo} still applies when the points $\{p_i\}$ are fixed under the $\cpx^*$-action; for example one can apply \ref{exa_proj} to fixed points with appropriate multiplicities. Moreover, we expect that \ref{main_teo} still holds for arbitrarily singular fixed points on a polarised variety $(X, L)$; of course $F(X)$ will need to be replaced with the Donaldson-Futaki invariant for the $\cpx^*$-action on the base. 
\end{rem}
\begin{rem}
A very interesting generalisation of \ref{main_teo} has been found by Della Vedova \cite{alberto}. Roughly speaking this applies to \textit{extremal} metrics, i.e. metrics whose scalar curvature has holomorphic $(1,0)$-gradient. See also \cite{arezzoIII} for more details. 
\end{rem}

\vskip.3cm

\noindent Universit\`a di Pavia, Via Ferrata 1 27100 Pavia, ITALY\\
and\\
Imperial College, London SW7 2AZ, UK.\\
\textit{E-mail:} jacopo.stoppa\texttt{@unipv.it}

\end{document}